\documentclass[12pt]{article}
\makeindex
\usepackage{epsfig,amsmath,a4wide,amssymb,amsfonts,amsbsy,xy,latexsym,epsf,pstricks,verbatim,times}
\xyoption{all}
\usepackage[english]{babel}

\newcommand{\PP}{\mathcal P}
\newcommand{\PQ}{\mathcal Q}

\let\set\mathbb
\def\<<{\leavevmode
  \raise0.28ex\hbox{$\scriptscriptstyle\langle\!\langle$}\nobreak
  \hskip -.6pt plus.3pt minus.2pt\,}
\def\>>{\,\nobreak\hskip -.6pt plus.3pt minus.2pt
  \raise0.28ex\hbox{$\scriptscriptstyle\rangle\!\rangle$}}

\def\vare{{\varepsilon}}
\def\dx{{d_x}}
\def\dy{{d_y}}

\def\AA{{\set A}}

\def\Hom{\mathop{\rm{Hom}}\nolimits }

\def\Ker{\mathop{\rm{Ker}}\nolimits }
\def\Gal{\mathop{\rm{Gal}}\nolimits }
\def\Pic{\mathop{\rm{Pic}}\nolimits }

\def\CC{{\set C}}

\def\End{\mathop{\rm End }}
\def\FF{{\set F}}

\def\Fp{{\FF _p}}
\def\Fps{{{\FF ^*_p}}}

\def\Fq{{\FF _q}}
\def\Fqs{{{\FF ^*_q}}}
\def\Fqd{{{\FF _{q^d}}}}

\def\bK{{\bf K  }}
\def\bL{{\bf L  }}

\def\bG{{\bf  G }}
\def\bH{{\bf H  }}
\def\oG{{\oplus_\bG}}
\def\oGm{{\oplus_{\bG_m}}}
\def\oGa{{\oplus_{\bG_a}}}
\def\oH{{\oplus_\bH}}
\def\zG{{0_\bG}}
\def\zH{{0_\bH}}

\def\PU{{{\set P}^1}}
\def\AA{{\set A}}
\def\AU{{{\set A}^1}}

\def\QQ{{\set Q}}

\def\ZZ{{\set Z}}
\def\agot{{\mathfrak a}}
\def\igot{{\mathfrak i}}

\def\bgot{{\mathfrak b}}
\def\cA{{\cal A}}
\def\cB{{\cal B}}
\def\cC{{\cal C}}

\def\cF{{\cal F}}
\def\cG{{\cal G}}
\def\cH{{\cal H}}
\def\cI{{\cal I}}

\def\cL{{\cal L}}

\def\cS{{\cal S}}

\def\cU{{\cal U}}

\def\cgot{{\mathfrak  c}}

\def\mmu{{\set \mu}}

\begin{document}
\author{Jean-Marc Couveignes\thanks{Institut de  Math\'ematiques  de Toulouse,
Universit\'e de Toulouse et CNRS.} and Reynald Lercier\thanks{Centre d'\'Electronique de l'Armement, 35170 Bruz, France.}}
\title{Galois invariant smoothness basis\thanks{Research supported by the
    French D{\'e}l{\'e}gation G{\'e}n{\'e}rale pour l'Armement,
Centre d'\'Electronique de l'Armement and by the Fonds National pour la
Science (ACI NIM).}}

\maketitle

\bibliographystyle{plain}

\begin{abstract}
This text answers a question raised by Joux and the second author  about
the computation of  discrete logarithms in the multiplicative group
of finite fields. Given a finite residue field  $\bK$, one looks for a
smoothness
basis for  $\bK^*$ that is left invariant by automorphisms of $\bK$.
For a broad class of  finite fields, we manage to construct
models that allow such a smoothness basis.
This work aims at accelerating discrete logarithm computations
in such fields. We treat the cases
of codimension  one (the linear sieve) and codimension two
(the function field sieve).
\end{abstract}

\hfill {\it To Gilles Lachaud, on the occasion of his 60th birthday}

%\tableofcontents

\section{Motivation}\label{section:presentation}

We look for finite fields that admit Galois invariant
smoothness basis. It is known that such basis accelerate the
calculation of discrete logarithms. We first
recall this observation by Joux and Lercier
in section~\ref{section:introduction} and we give a
first example of this situation in  section~\ref{section:exemple}. 
We recall in section~\ref{section:KAS} the
rudiments of Kummer and Artin-Schreier theories. These theories
produce the known examples of such smoothness basis.
We then show in  section~\ref{section:espaces} that the only
extensions admitting Galois invariant flags of linear spaces
are given by those two theories. In section
~\ref{section:isog}, we consider   a more general
setting: specialization of isogenies
between algebraic groups. We deduce a first
non trivial example of Galois invariant smoothness basis in
section~\ref{section:tore}. 
In the next section~\ref{section:ell}, we show that elliptic
curves produce a range of such invariant basis, provided the
degree of the field is not too large. 

In section~\ref{section:JL}, we recall
the principles of fast sieving algorithms (the number field sieve and
the function field sieve). We 
show in 
section
\ref{section:carre} that our approach can be adapted to these algorithms.
A detailed example is given in section~\ref{section:exp}.
We finish with a few remarks and questions about the relevance  of our method.

\section{A remark by Joux and Lercier}\label{section:introduction}

We recall in this section the principle
of a simple algorithm
for computing discrete logarithms in the multiplicative
group of a finite field $\Fq$ where  $q=p^d$ and $d\ge 2$.
 See~\cite{odlyzko} for a survey
on discrete logarithm computation.

The finite field  $\Fq$  is seen as a residue
field $\Fp[X]/(A(X))$ where   $A(X)\in \Fp[X]$ is
a degree $d$ unitary irreducible polynomial. We set
$x=X\bmod A(X)$. 
Let  $k$ be an integer such that  $0\le k\le d-1$ and 
let  $V_k\subset \Fq$ be the $\Fp$-vector space generated
by  $1$, $x$, \dots, $x^{k}$. 
So  $V_0=\Fp\subset V_1\subset \ldots \subset V_{d-1}=\Fq$ and  
$V_k\times V_l\subset V_{k+l}$ if  $k+l\le d-1$.

One looks  for multiplicative relations between 
elements of  $V_\kappa$ for some integer  $\kappa$. 
For example, if one  takes  $\kappa=1$, the  relations
we are looking  for take the form 
\begin{equation}\label{eqn:prod}
\prod_{i} (a_i+b_ix)^{e_i}=1\in \Fq
\end{equation}
where  the $a_i$ and  $b_i$ lie in  $\Fp$.
We collect such relations until
we obtain a basis
of the 
$\ZZ$-module of relations between
elements in  $V_\kappa$.

How do we find relations like relation~(\ref{eqn:prod}) ?
Assume again  $\kappa=1$.
The simplest form of the sieving algorithm
picks  random triplets $(a_i,b_i,e_i)$ and computes
the remainder  $r(X)$ of the Euclidean
division of  $\prod_i(a_i+b_iX)^{e_i}$ by  $A(X)$. So 
$$r(X)\equiv \prod_i(a_i+b_iX)^{e_i}\bmod A(X)$$
where $r(X)$ is a more or less random
polynomial in $\Fp[X]$ with degree $\le  d-1$.

We hope  $r(X)$  decomposes as a product of polynomials
with degree smaller than or equal to  $\kappa=1$. 
If this is the case, we find  $r(X)=\prod_j(a'_j+b'_jX)^{e'_j}$ and we obtain
a relation 
$$\prod_{ i } (a_i+b_ix)^{e_i}\prod_j(a'_j+b'_jx)^{-e'_j}=1$$
of the expected form.
One says  that $V_\kappa$ is the smoothness basis.

Joux and  Lercier  notice in~\cite{JouxLercier2} that,
if there exists an  automorphism 
$\agot$ of  $\Fq$ such that  $\agot (x)= ux+v$ with  $u$, $v \in \Fp$,
then the action of  $\agot$ on equation~(\ref{eqn:prod}) produces
another equation of the same kind. 
Since the efficiency of
discrete logarithm algorithms depends 
on the number of such equations one can produce
in a given amount of time, one wishes to know when
such useful automorphisms exist.
We also wonder how to generalize 
this  observation.

We stress that  $\agot$ acts
both on equations and factors of the form  $a_i+b_ix$. Rather than
increasing the number of equations, such an action may be used
to lower the number of factors involved in them.
If  $\agot$ is the  $n$-th power
of  the Frobenius automorphism, we obtain for free 
\begin{displaymath}
\agot(a+bx)=(a+bx)^{p^n}=v+a+ubx
\end{displaymath}
So we can remove  $v+a+ubx$ out of the smoothness basis
and replace it everywhere by  $(a+bx)^{p^n}$.
This way, we only keep one element in every 
orbit of
the Galois group acting on 
$V_\kappa$.
As a consequence, the size of the
linear system we must solve is divided by the order
of the group generated by  $\agot$. If  $\agot$
generates the full  Galois group of  $\Fq/\Fp$, then the number
of unknowns is divided by $d$, 
the degree of the finite field  $\Fq$.

Our concern in
this text is to find 
models for finite fields for which the automorphisms 
respect the special  form of   certain elements.
For example, if the finite field
is given as above, the elements are given as polynomials
in 
 $x$.   Any element  $z$ of the finite field
has a degree: This is the smallest integer 
 $k$ such that  $z\in V_k$. The degree 
of   $a_0+a_1x+\dots +a_kx^k$ is thus  $k$ provided
 $0\le k<d$ and  $a_k\not =
0$ (and by convention, $\deg 0 = 0$ ).
The degree is sub-additive,
 $\deg(w\times z)\le \deg(w)+\deg(z)$.

The question raised boils down to
asking if this degree function is preserved 
by the automorphisms of  $\Fq$.
It is worth noticing that
the interest of the degree function in this context
comes from the following properties.
\begin{itemize}
\item  The degree is sub-additive
(and often  even  additive): The degree of the
product of two non zero elements is the sum of the degrees
of either elements 
provided this sum is  $<d$.
\item The degree sorts nicely  the
elements of $\Fq$: 
There are   $q^n$ elements of degree 
 $<n$ if  $1\le n\le d$.
\item There exists a factoring  algorithm that 
decomposes some elements in  $\Fq$ as products
of elements with smaller degrees (\textit{e.g.} with degree $\le \kappa$). The density
of such $\kappa$-smooth elements is not too small.
\end{itemize}

In this article, we look for
such degree functions on finite fields
having the extra property that they are Galois invariant: Two conjugate
elements have the same degree.

\section{A first example}\label{section:exemple}

Here is an example provided by
 Joux and   Lercier.
Take  $p=43$ and  $d=6$, so  $q=43^6$, and set  $A(X)=X^6-3$ which is
an irreducible polynomial in  $\FF_{43}[X]$. So  $\Fq$ is
seen as the residue field
$\FF_{43}[X]/(X^6-3)$. 

One checks that 
$p=43$ is congruent to  $1$ modulo
$d=6$, so 
$\phi(x)= x^{43} = (x^6)^7\times x=3^7x=\zeta_6 x$
where  $\zeta_6=3^7=37\bmod 43$ is
a primitive sixth root of unity.
Since the  Frobenius $\phi$ generates the Galois
group, one can divide by $6$ the size
of the smoothness basis.

In the second example
provided by  Joux and  Lercier (and coming from XTR of type T30) one
takes $p=370801$ and  $d=30$ with $A(X)=X^{30}-17$.
This time,  $p$ is congruent to
$1$ modulo $d=30$  and 
$\phi(x)=x^p=x^{30\times 12360}\times x=\zeta_{30}x$
where $\zeta_{30}=17^{12360}\bmod p=172960\bmod p$.
As a consequence, one can
divide by $30$ the size of the smoothness basis.

We are here in the context
of Kummer theory. In the next section we
recall the basics of this theory, that classifies
cyclic extensions of $\Fp$ with degree $d$ dividing $p-1$. 
Artin-Schreier theory is the counterpart for  cyclic
 $p$-extensions in characteristic $p$ and
we sketch it as well.
Both theories are of very limited interest for our purpose.
We shall need to consider 
the more general situation of an algebraic group with
rational torsion.

\section{ Kummer and  Artin-Schreier theories}\label{section:KAS}

The purpose here is
to classify cyclic extensions
of degree $d\ge 2$ of a field $\bK$ with characteristic $p$ 
in two simple cases.
\begin{itemize}
\item Kummer case: 
$p$ is prime to  $d$ and $\bK$ contains a primitive
$d$-th root of unity;
\item Artin-Schreier case:  $d=p$.
\end{itemize}

\paragraph*{Kummer theory.}
We follow Bourbaki~\cite[A V.84]{bourbaki}.
According to  Kummer theory, if  $p$ is prime to  $d$ 
and  $\bK$ contains  a primitive $d$-th root of unity,
then every degree $d$ cyclic extension of $\bK$
is generated by a radical.

Assume $\bK$ is embedded in some algebraic closure $\bar \bK$.
To every $a$ in  $\bK^*/(\bK^*)^d$ (which we may
regard as an element in  $\bK^*$), we associate the field 
 $\bL=\bK(a^\frac{1}{d})$ where
 $a^\frac{1}{d}$ is any root of
$X^d-a$ in  $\bar\bK$.

The map  $x\mapsto x^d$ is an epimorphism
from the multiplicative group $\bar\bK^*$ onto itself.
The kernel of this epimorphism
is the group of $d$-th roots of unity.
The roots of $X^d-a$ lie in the inverse image of
$d$ by this epimorphism.

The field  $\bK(a^\frac{1}{d})$ may not
be isomorphic to $\bK[X]/(X^d-a)$. It is when
 $a$ has order $d$ in the group
$\bK^*/(\bK^*)^d$.
On the other hand, if  $a$ lies in  $(\bK^*)^d$ then 
 $\bK[X]/(X^d-a)$   is the product
of  $d$ copies of  $\bK$.

Let's come back to the case
when  $a$ has order  $d$ in $\bK^*/(\bK^*)^d$. 
The degree $d$ extension  $\bL/\bK$ is Galois since, if
we set  $b=a^\frac{1}{d}$, we have 
$$X^d-a=(X-b)(X-b\zeta_{d})(X-b\zeta_{d}^2)\dots(X-b\zeta_{d}^{d-1})$$
where  $\zeta_d$ is a primitive
$d$-th root of unity.
The Galois group of $\bL/\bK$ 
is made of transformations 
$\agot_n: x\mapsto x\zeta_d^n$
and the map
$n\mapsto \agot_n$ is an isomorphism from the group $\ZZ/d\ZZ$ onto
$\Gal(\bL/\bK)$.

To avoid distinguishing too many cases,
one follows Bourbaki~\cite[A V.84]{bourbaki}.
Rather than a single
element in  $\bK^*/(\bK^*)^d$ one picks
a subgroup  $H$ of  $\bK^*$ containing  $(\bK^*)^d$ and
one forms the extension  $\bK(H^\frac{1}{d})$ by adding
to $\bK$ all $d$-th roots of all elements in $H$.
To every automorphism  $\agot$ in  $\Gal(\bK(H^\frac{1}{d})/\bK)$, 
one associates an homomorphism 
 $\psi(\agot)$  from $H/(\bK^*)^d$ to the group  $\mmu_d(\bK)$
of $d$-th roots of unity. The homomorphism 
$\psi (\agot)$  is defined by 
$$\psi(\agot): \theta \mapsto
\frac{\agot(\theta^\frac{1}{d})}{\theta^\frac{1}{d}}$$
where  $\theta^\frac{1}{d}$ is one of
the $d$-th roots of $\theta$.
The map  $\agot \mapsto \psi (\agot)$ is an isomorphism from the
$\Gal(\bK(H^\frac{1}{d})/\bK)$
onto $\Hom(H/(\bK^*)^d, \mmu_d(\bK))$.
This presentation
of Kummer theory constructs 
abelian extensions
of  $\bK$ with exponent dividing  $d$.

In the case we are interested in, the field
$\bK=\Fq$ is finite.
Any subgroup  $H$ of  $\bK^*$ is cyclic. 
In order to have $\mmu_d$ in $\bK$, one assumes
that $d$ divides  $q-1$.
We set  $q-1=md$. The group $(\bK^*)^d$ has order  $m$. 
The  quotient $\bK^*/(\bK^*)^d$ is cyclic
of order   $d$. It is natural to take  $H=\bK^*$.
We find  the unique degree $d$  cyclic extension  $\bL$
of $\bK$. 
It is generated by a $d$-th root of a generator
$a$ of $\bK^*$.

Set  $b=a^\frac{1}{d}$ and  $\bL=\bK(b)$. 
The Galois group  $\Gal(\bL/\bK)$ is generated
by the  Frobenius
$\phi$ and the action of  $\phi$ on  $b$ is given by 
$\phi(b)=b^{q}$, so 
$$\zeta=\frac{\phi(b)}{b}=b^{q-1}=a^m$$
is a  $d$-th  root of
unity that depends on  $a$. The map  $a\mapsto \zeta$ 
is an  isomorphism of  $\bK^*/(\bK^*)^d$  onto
 $\mmu_d(\bK)$ which is nothing but exponentiation by $m$.

The limitations of
this construction are
clear: It requires
primitive $d$-th roots of unity in  $\bK$.
Otherwise, one may
jump to some auxiliary extension 
$\bK'=\bK(\zeta_d)$  of  $\bK$,  that
may be quite large. One applies 
Kummer theory to this bigger extension and one obtains 
a degree $d$ cyclic extension $\bL'/\bK'$.
Descent can be performed using
resolvants (see~\cite[Chapter III.4]{matzat}) at a serious
computational expense. We shall not follow this track.

\subparagraph*{Example.}
Coming back to the first
example one finds  $q=p=43$, $p-1=42$, $d=6$, $m=7$,
$a=3$ and  ${\phi(b)}/{b}=a^m=3^7\bmod 43$.

\paragraph*{Artin-Schreier theory.}

We follow 
Bourbaki~\cite[A V.88]{bourbaki}. 
If $p$  is the characteristic of  $\bK$,
then any cyclic degree $p$  extension of $\bK$ is
generated by the roots of a polynomial of the
form
$$X^p-X-a=\wp(X)-a=0$$
where  $a\in \bK$ and the expression $\wp(X)=X^p-X$ 
plays a similar role to  $X^d$ in Kummer theory.
The map  $x\mapsto \wp(x)$ defines an epimorphism
from the additive group  $\bar\bK$ onto itself. 
The kernel of this epimorphism
is the additive group of
the prime field  $\Fp\subset \bar\bK$.

Let  $a$ be an element of  $\bK/\wp(\bK)$ (that we may see
as an element of   $\bK$ in this class). One associates
to it the extension field  $\bL=\bK(b)$ where $b\in\wp^{-1}(a)$. 
If $a$ has order $p$ in  $\bK/\wp(\bK)$, the 
extension $\bL/\bK$ has degree  $p$  and is Galois since we have 
$$X^p-X-a=(X-b)(X-b-1)(X-b-2)\dots(X-b-(p-1)).$$
The  Galois group is made of transformations
of the form 
$\agot_n: x\mapsto x+n$
and the map
$n\mapsto \agot_n$ is an  isomorphism from the group  $\ZZ/p\ZZ$ 
onto 
$\Gal(\bL/\bK)$.

Again, if one wishes to construct
all abelian extensions of $\bK$
with  exponent $p$ one follows  Bourbaki~\cite[A V.88]{bourbaki}. 
One takes a subgroup  $H$ of  $(\bK,+)$ containing $\wp(\bK)$ 
and one forms the extension  $\bK(\wp^{-1}(H))$.
To every automorphism  $\agot$ in $\Gal(\bK(\wp^{-1}(H))/\bK)$,
one associates an  homomorphism $\psi(\agot)$ from  $H/\wp(\bK)$ 
onto the additive  group  $\Fp$ of the prime field. The homomorphism
$\psi (\agot)$  is defined by 
$$\psi(\agot): \theta \mapsto
\agot(c)- c$$
where  $c$ belongs to $\wp^{-1}(\theta)$,
the fiber of  $\wp$ above  $\theta$.

The map  $\agot \mapsto \psi (\agot)$ is an isomorphism 
from the  Galois  group 
$\Gal(\bK(\wp^{-1}(H))/\bK)$
onto  $\Hom(H/\wp(\bK), \Fp)$.

In our case, the field  $\bK=\Fq$ is finite of characteristic $p$.  
We set  $q=p^f$.
The morphism 
$\wp: \Fq \rightarrow \Fq$
has kernel
$\Fp$ and the quotient  $\Fq/\wp(\Fq)$ has order  $p$.
The unique
degree $p$ extension  $\bL$ of  $\Fq$ 
is generated by  $b\in \wp^{-1}(a)$ where  $a\in \Fq-\wp(\Fq)$.
The   Galois group $\Gal(\bL/\bK)$ is generated by the
 Frobenius $\phi$  and  
$\phi(b)-b$ belongs to  $\Fp$.  The map $a\mapsto \phi(b)-b$ 
is an  isomorphism from 
 $\bK/\wp(\bK)$  onto   $\Fp$.

Let us make this isomorphism more explicit. 
We have   $\phi(b)=b^q$ where  $q=p^f$ is the order of  $\bK=\Fq$. One computes
$$\phi(b)-b=b^q-b=(b^p)^{p^{f-1}}-b=(b+a)^{p^{f-1}}-b\text{ since }\wp(b)=b^p-b=a.$$
So  $b^{p^{f}}-b=b^{p^{f-1}}-b +a^{p^{f-1}}$. Iterating, we 
obtain  
$$\phi(b)-b=b^{p^{f}}-b=a+a^p+a^{p^2}+\cdots+a^{p^{f-1}}.$$
The isomorphism from  $\bK/\wp(\bK)$  onto the additive group   $\Fp$
is nothing but the absolute trace.

\subparagraph*{Example.} Take  $p=7$ and  $f=1$,
so $q=7$. The absolute  trace of $1$ is  $1$,
so we set  $\bK=\FF_7$ and  $A(X)=X^7-X-1$ and we set 
$\bL=\FF_{7^7}=\FF_7[X]/(A(X))$.  Setting  $x=X\bmod A(X)$, one has
 $\phi(x)=x+1$.

\section{Invariant linear
spaces of a cyclic extension}\label{section:espaces}

Let us recall that the question raised in
section~\ref{section:introduction} concerns
the existence of automorphisms
that stabilize a given smoothness basis.
We saw that smoothness basis
are usually made using flags of linear spaces.
Therefore, one wonders if, for a given cyclic extension
 $\bL/\bK$, there exists
 $\bK$-vector subspaces of  $\bL$ that
are left invariant by the Galois group of  $\bL/\bK$.

Let   $d\ge 2$ be  an integer 
and $\bL=\bK[X]/(X^d-r)$  a Kummer extension. For any integer
 $k$ between  $0$ and  $d-1$, let 
$L_k=\bK\oplus \bK x\oplus \cdots\oplus \bK x^k$
be the  $\bK$-vector subspace generated by
the $k+1$ first powers of  $x=X\bmod X^d-r$.
The  $L_k$ are invariant under Galois action since for  $\agot$,
a 
$\bK$-automorphism of  $\bL$,
there exists a  $d$-th root of unity  $\zeta \in \bK$
such that 
$$\agot(x)=\zeta x$$
 and   $\agot(x^k)=\zeta^k x^k$.
One has a flag of  $\bK$-vector spaces,
$V_0=\bK\subset V_1\subset \dots \subset V_{d-1}=\bL$,
respected by Galois action. 
So the  ``degree'' function is invariant
under this  action.
This is exactly what happens
in the two examples of  section~\ref{section:introduction}.
If the smoothness basis is made of 
irreducible polynomials  of degree 
 $\le \kappa$, then it is acted on by the Galois group.

If now  $\bL=\bK[X]/(X^p-X-a)$ is an Artin-Schreier
extension, for every integer
 $k$ between  $0$ and  $p-1$, we call
$V_k=\bK\oplus \bK x\oplus \cdots\oplus \bK x^k$
 the  $\bK$-vector space generated by
the $k+1$ first powers
of  $x=X\bmod X^p-X-a$.
The  $V_k$ are  globally invariant 
under Galois action. Indeed, if $\agot$ is a 
$\bK$-automorphism of $\bL$,  then
there is a    $n \in \Fp$
such that
$\agot(x)=x+n,$
 so
$$\agot(x^k)=(x+n)^k=\sum_{0\le \ell\le k}\left( \begin{smallmatrix}   
k\\ \ell
\end{smallmatrix}  \right)n^{k-\ell}x^\ell.$$

We find again a flag of
 $\bK$-vector spaces,
$V_0=\bK\subset V_1\subset \dots \subset V_{p-1}=\bL$,
that is fixed by Galois action.
This time, the Galois action
is no longer diagonal
but triangular.
For cyclic extensions of degree
a power of  $p$,  Witt-Artin-Schreier theory also
produces a flag of Galois invariant vector spaces.
See the beginning  of  Lara Thomas's thesis
\cite{lara}  for an introduction with references.
\medskip

One may wonder if Galois invariant
flags of vector spaces
exist for other  cyclic field extensions.
Assume  $\bL/\bK$ is a degree $d$ cyclic
extension where $d$ is prime to the characteristic
$p$.
Let $\phi$ be a generator of
the Galois group  $C=<\phi>=\Gal(\bL/\bK)$.
 According to the normal
basis theorem~\cite[Theorem 13.1.]{langalgebra}, there
exists a  $w$  in  $\bL$ such that
$$(w,\phi(w),\phi^2(w),\dots,\phi^{d-1}(w))$$
 is a  $\bK$-basis of $\bL$.
Therefore $\bL$, as a $\bK[C]$-module, is isomorphic
to the regular representation. 
The order  $d$ of $C$ being prime to the characteristic,
the ring  $\bK[C]$ is semi-simple according to
 Maschke theorem~\cite[Theorem 1.2.]{langalgebra}. 
The characteristic polynomial of
 $\phi$ acting on the $\bK$-vector space $\bL$ is
$X^d-1$. This is a separable polynomial in  $\bK[X]$.

To every  $\bK$-irreducible factor $f(X)\in \bK[X]$ of  $X^n-1$,
there corresponds a  unique irreducible characteristic
subspace $V_f\subset \bL$,  invariant by  $\phi$. The
characteristic polynomial of  $\phi$ restricted to  $V_f$ 
is $f$.
According to
Schur's lemma~\cite[Proposition 1.1.]{langalgebra},
any $\bK[C]$-submodule of
$\bL$ is a direct sum of 
some $V_f$.

Assume there exists a complete flag
of $\bK$-vector spaces, each invariant
by $\phi$,
$V_0=\bK\subset V_1\subset \dots \subset V_{d-1}=\bL$,
where $V_k$ has  dimension $k$. Then
all irreducible factors of $X^d-1$ must have
degree $1$.
So $\bK$ contains primitive roots of unity
and we are in the context of  Kummer theory.
To every Galois invariant flag, there corresponds an
order on $d$-th roots of unity (or equivalently
on the associated characteristic spaces in $\bL$). 
There are  $d!$ such flags. 

The flags produced by Kummer theory are of the following form:
\begin{multline*}
  V_1\subset V_1\oplus V_\zeta\subset V_1\oplus V_\zeta\oplus V_{\zeta^2}\subset  \dots\\ \subset 
V_1\oplus V_\zeta\oplus V_{\zeta^2}\oplus\cdots\oplus V_{\zeta^{d-2}}
\subset V_1\oplus V_\zeta\oplus V_{\zeta^2}\oplus\cdots\oplus V_{\zeta^{d-2}}\oplus
V_{\zeta^{d-1}}
\end{multline*}
where $\zeta $ is a  primitive $d$-th root of
unity
and  $V_\zeta$ is  $V_{X-\zeta}$, the eigenspace
associated to $\zeta$.

Among the $d!$ flags that are  $\phi$-invariants, only
$\varphi (d)$ come from Kummer theory.
They correspond to the  $\varphi(d)$ primitive roots
of unity. These flags enjoy a multiplicative property:
If  $k\ge 0$ and  $l\ge 0$
and  $k+l\le d-1$, then
$V_k\times V_l\subset V_{k+l}.$

The conclusion of this section
is thus rather negative. If we want to go further than
Kummer theory, we cannot ask for Galois invariant flags
of vector subspaces.

\section{Specializing isogenies between
commutative algebraic groups}\label{section:isog}

Kummer and Artin-Schreier
theories are two special cases of a general situation that
we now describe. Our aim is to produce nice models
for a broader variety of finite fields.

Let  $\bK$ be a field and  $\bG$  a commutative  algebraic group
over $\bK$.
Let  $T\subset \bG(\bK)$ be a non trivial finite group
of $\bK$-rational points in  $\bG$ and let 
$$I: \bG\rightarrow \bH$$
be the quotient isogeny of  $\bG$ by  $T$.
Let  $d\ge 2$ be the  cardinality of  $T$ which is also
the degree of $I$.
Assume there exists a  $\bK$-rational point  $a$ on 
$\bH$ such that  $I^{-1}(a)$ is irreducible
over  $\bK$. Then every point  $b\in \bG(\bar \bK)$
such that  $I(b)=a$ defines a cyclic degree $d$
 extension $\bL$ of  $\bK$: We set  $\bL=\bK(b)$
and we notice that the  geometric
origin of this extension results in a nice
description of   $\bK$-automorphisms  of $\bL$.

Let  $t$ be a point  in  $T$ and let  $\oG$ stand 
for the addition law in the algebraic 
 group   $\bG$.  Let  $\oH$ stand for the addition law in $\bH$.
We denote by  $\zG$ the unit element
in  $\bG$ and  $\zH$ the one in 
$\bH$.
The point  $t\oG b$ verifies 
$$I(t\oG b)=I(t)\oH I(b)=\zH\oH a=a.$$

So  $t\oG b$ is Galois conjugated to  $b$ and all conjugates
are obtained that way from all  points  $t$ in   $T$.
So we have an  isomorphism
between $T$ and  $\Gal(\bL/\bK)$, which associates
to every  $t\in T$  the residual automorphism 
$$b\in I^{-1}(a)\mapsto b\oG t.$$

Now,  assuming the geometric 
formulae that describe the translation 
$P\mapsto P\oG t$ in  $\bG$ are  simple enough, we obtain
a nice description 
of the  Galois group of  $\bL$ over $\bK$.

Kummer and Artin-Schreier theories 
provide two illustrations of this general geometric situation.

The algebraic group underlying 
 Kummer  theory is
the multiplicative group $\bG_m$ over the base field $\bK$. The isogeny  $I$ is the multiplication
by  $d$: 
$$I=[d]: \bG_m\rightarrow \bG_m.$$

One can see the group $\bG_m$ as a sub-variety
of the affine line   $\AA^1$ with $z$-coordinate. The
inequality $z\not =0$ defines the open subset $\bG\subset \AA^1$.
The origin  $\zG$ has  coordinate  $z(\zG)=1$. The  group law is given by
$$z(P_1\oGm P_2)=z(P_1)\times z(P_2).$$

Here we have  $\bH=\bG=\bG_m$ and the isogeny  $I$ can be given
in terms of the $z$-coordinates by 
$$z(I(P))=z(P)^d.$$

Points in the kernel of  $I$ have  for $z$-coordinates  the $d$-th roots
of unity.
The inverse image by  $I$ of a point $P$ in $\bG$ is made
of  $d$ geometric points having 
for $z$-coordinates the  $d$-th roots of  $z(P)$.
Translation by an element  $t$ of the kernel of  $I$, 
$P\mapsto P\oGm t$,
can be expressed in terms of  $z$-coordinates by 
$$z(P\oGm t)=z(P)\times \zeta$$
where  $\zeta=z(t)$ is the  $d$-th root
of unity associated by $z$ to the $d$-torsion point  $t$.

As far as  Artin-Schreier theory is concerned, the underlying algebraic  group 
is the additive group $\bG_a$ over the base
field $\bK$, identified with the affine 
line   $\AA^1$ over $\bK$. A point $P$ on  $\bG_a$ is given
by its $z$-coordinate. 
The origin  $\zG$ has coordinate  $z(\zG)=0$ and the group
law is given by 
$$z(P_1\oGa P_2)=z(P_1)+ z(P_2).$$

The degree $p$ isogeny  $I$ is
$\wp: \bG_a\rightarrow \bG_a$,
given in terms of $z$-coordinates by
$$z(\wp(P))=z(P)^p-z(P).$$

Here again  $\bH=\bG$. The $z$-coordinates
of  points in the kernel of
 $\wp$ are the elements of the prime field $\Fp$.
The inverse image by  $I$ of a point $P$ in  $\bG$ is made
of  $p$ geometric points whose $z$-coordinates are the
 $p$ roots of the equation $X^p-X=z(P)$.
Translation by an element  $t$ in the kernel of  $I$,
$P\mapsto P\oGa t$,
can be expressed in terms
of  $z$-coordinates by 
$$z(P\oGa t)=z(P) + \tau
\text{ where }\tau=z(t)\in \Fp\ .$$

\section{A different example}\label{section:tore}

We plan to apply the generalities in the previous section
to various algebraic groups. We guess every commutative algebraic group
may bring its contribution to the construction of Galois
invariant smoothness basis.  Since we look for  simple  
translation formulae,
we expect the simplest algebraic groups to be the most useful.
We start with the most familiar algebraic groups (after $\bG_m$ and $\bG_a$):
These are the dimension $1$ tori.
Let  $\bK$ be a field with characteristic different from  $2$ and let 
$D$ be a non zero element in  $\bK$.
Let  $\PU$ be the projective line with projective coordinates
 $[U,V]$. Let $u=\frac{U}{V}$ be the associated affine coordinate.
We denote by  $\bG$ the open subset of  $\PU$ defined by the inequality 
$$U^2-DV^2\not =0.$$

To every point  $P$ of  $\bG$, we associate its  $u$-coordinate, possibly
infinite
but distinct from  $\sqrt D$ and  $-\sqrt D$.
The unit element in  $\bG$ is the  point $\zG$ with projective 
coordinates $[1,0]$ and  $u$-coordinate $\infty$.
For $P_1\not = \zG$ and $P_2\not = \zG$, the addition law is given by 
$$u(P_1\oG P_2)=\frac{u(P_1)u(P_2)+D}{u(P_1)+u(P_2)}
\text{ and }
u(\ominus_\bG P_1)=-u(P_1).$$

We now assume  that $\bK=\Fq$ is a finite field and $D\in \Fqs$ is not a square in
 $\Fq$.
The group $\bG(\Fq)$ has order $q+1$ and the corresponding values
of  $u$ lie in  
$\Fq\cup \{\infty \}$.
The  Frobenius endomorphism, $\phi: \bG \rightarrow \bG$,
\begin{math}
[U,V] \rightarrow [U^q,V^q],
\end{math}
is nothing but multiplication by
$-q$. Indeed, let  $P$ be a  point with projective coordinates  $[U,V]$. The
projective coordinates of $R=[q]P$ are the  coordinates in  $(1, \sqrt{D})$ of 
$$(U+ V\sqrt{D})^q=U^q-\sqrt{D}V^q$$
because  $D$ is not a square in  $\Fq$.
So  $R$ has coordinates $[U^q,-V^q]$ and it is the inverse of $\phi(P)$.

We pick an integer $d\ge 2$ such that the $d$-torsion $\bG[d]$ is
$\Fq$-rational. This is equivalent to the condition that $d$ divides $q+1$.
We set $q+1=md$.
Let  $I$ be the multiplication by 
$d$ isogeny,
$I=[d]: \bG\rightarrow \bG$,
with kernel  the cyclic group $\bG[d]$ of order $d$.
The quotient  $\bG(\Fq)/I(\bG(\Fq))=\bG(\Fq)/\bG(\Fq)^d$
is cyclic of order  $d$.

Let $a$ be a generator of  $\bG(\Fq)$ and $b$ a geometric point
in the fiber of $I$ above $a$. Let 
$u(b)$  be the  $u$-coordinate  of  $b$ and set 
$\bL=\bK(u(b))$. This is a degree  $d$ extension of  $\bK=\Fq$.
So  $\bL=\Fqd$.

The  Galois group of  $\Fqd/\Fq$  is isomorphic to  
$\bG[d]$:  For any  $\agot \in \Gal (\Fqd/\Fq)$, the  difference
$\agot(b)\ominus_\bG b$ is in  $\bG[d]$ and the pairing
$$(\agot , a)\mapsto  \agot(b)\ominus_\bG b$$
defines an isomorphism of  
$\Gal(\Fqd/\Fq)$ onto $\Hom(\bG(\Fq)/(\bG(\Fq))^d,\bG[d])$.

Here  $\Gal(\Fqd/\Fq)$ is cyclic of order $d$ generated by the  Frobenius
$\phi$. The pairing 
$(\phi,a)$ equals  $\phi(b)\ominus_\bG b$.
Remember that  $\phi(b)=[{-q}]b$  in $\bG$. So 
\begin{equation}\label{eqn:frob}
(\phi,a)=[-q-1]b=[-m]a.
\end{equation}

We obtain an exact description of Galois action on 
$I^{-1}(a)$.
It is given by   translations of the form $P\mapsto P\oG t$ with   $t\in
\bG[d]$. If we denote by $\tau$ the affine coordinate of $t$ and by $u$ 
the coordinate of  $P$ then the action is given by 
$$u\mapsto \frac{\tau u+D}{u+\tau},$$
which is rather nice since it is a rational linear transform.

We form the polynomial 
$$A(X)=\prod_{b\in I^{-1}(a)} (X-u(b))$$
 annihilating the  $u$-coordinates
of points in the inverse image of  $a$ by $I$.
This is a  degree $d$ polynomial
with coefficients in $\bK=\Fq$.
It is irreducible in $\Fq[X]$ because  $a$ generates  $\bG (\Fq)$. 
We have  $\bL=\bK[X]/(A(X))=\Fqd$.

The exponentiation formulae
in  $\bG$ give the explicit form of  $A(X)$.
One has
$$(U+\sqrt DV)^d=\sum_{0\le 2k\le d} 
\left(\begin{smallmatrix}
d\\2k \end{smallmatrix}\right)
U^{d-2k}V^{2k}D^k+\sqrt D\sum_{1\le 2k+1\le
  d} 
\left(\begin{smallmatrix}
d\\2k+1 \end{smallmatrix}\right)
U^{d-2k-1}V^{2k+1}D^k.$$
So,
$$u([d]P)=\frac{\sum_{0\le 2k\le d} u(P)^{d-2k}
\left(\begin{smallmatrix}
d\\2k \end{smallmatrix}\right)
D^k}{\sum_{1\le 2k+1\le
  d} u(P)^{d-2k-1}
\left(\begin{smallmatrix}
d\\2k+1 \end{smallmatrix}\right)
D^k}.$$
And 
$$A(X)=\sum_{0\le 2k\le d} X^{d-2k}
\left(\begin{smallmatrix}
d\\2k \end{smallmatrix}\right)
D^k-u(a)\sum_{1\le 2k+1\le
  d} X^{d-2k-1}
\left(\begin{smallmatrix}
d\\2k+1 \end{smallmatrix}\right)
D^k.$$

We set  $x=X\bmod A(X)$. Since
every $\Fq$-automorphism of $\Fqd$ transforms
$x$ into a linear rational fraction of
 $x$, it is natural to define for every
integer $k$
such that $k\ge 0$ and   $k< d$ the subset 
$$V_k=\{ \frac{u_0+u_1x+u_2x^2+\cdots+u_kx^k}{v_0+v_1x+v_2x^2+\cdots+v_kx^k}
  \ |\ (u_0,u_1,\ldots,u_k,v_0,v_1,\ldots,v_k)\in \bK^{2k+2}  \}\ .$$

One has
$\Fq=V_0\subset V_1\subset \dots\subset V_{d-1}=\Fqd$
and the  $V_k$ are  Galois invariant.
Further, it is clear that 
$V_k\times V_l\subset V_{k+l}$
provided  $k+l\le d-1$.
Again we find a flag of
Galois invariant subsets
of $\bL=\Fqd$. But these subsets are no longer
vector spaces.

If we define the degree
of an element of  $\bL$ to be the smallest integer 
 $k$ such that  $V_k$ contains this element, then the degree
is Galois invariant and sub-additive,
$\deg(wz)\le \deg(w)+\deg(z).$
The degree this times takes values between  $0$ and  $\lceil \frac{d-1}{2}
\rceil $. It is a slightly less informative function than in the  Kummer 
or Artin-Schreier cases (it takes twice less values).

\subparagraph*{Example.} Take  $p=q=13$ and  $d=7$. So  $m=2$.
Let   $D=2$ which is not  a square in 
$\FF_{13}$.
We look for some $a=U+\sqrt 2 V$ such that  $U^2-2V^2=1$ and 
$a$ has order $p+1=14$ in $\FF_{13}(\sqrt 2)^*$.
For example  $U=3$ and  $V=2$ are fine. The
$u$-coordinate  of  $3+2\sqrt 2$ is  $u(a)=\frac{3}{2}=8$.
One can write the polynomial
$$A(X)=X^7+3X^5+10X^3+4X-8(7X^6+5X^4+6X^2+8).$$
Formula~(\ref{eqn:frob}) predicts the Frobenius action.
We set  $t=[-m]a=[-2]a$ so  $u(t)=4$ and  Frobenius operates by translation by 
$t$, so
$X^p= \frac{4X+2}{X+4} \bmod A(X).$
\bigskip

So we have made a small progress: We can now treat extensions of $\Fq$ of
degree
dividing $q+1$. Unfortunately this condition is just as restrictive (though
different) as the one
imposed by Kummer theory. What do we do if  the degree does not divide
$q+1$ nor  $q-1$ ?

We must diversify the algebraic groups we use. The next family to consider is
made of elliptic curves.

\section{Residue fields of divisors on elliptic curves}\label{section:ell}

We now specialize the computations in  section~\ref{section:isog} to the case
where  $\bG$ is an elliptic curve.
Take  $\bK=\Fq$ a finite field for which we want to construct a degree 
$d\ge 2$
extension where  $d$ is prime to the characteristic  $p$ of $\Fq$. 
Here  $\bG=E$ is an ordinary elliptic curve over  $\Fq$. 
We denote by $\phi$ the Frobenius endomorphism of $E$.
Let  $\igot$ be an
invertible ideal in the endomorphism ring $\End(E)$ of $E$. Assume $\igot$
divides   $\phi-1$ and $\End(E)/\igot$ is cyclic of order  $d\ge 2$.
So $E(\Fq)$ contains
a  cyclic subgroup  $T=\Ker \igot$ of order $d$.

Let  $I: E \rightarrow F$ be the degree $d$ cyclic isogeny with kernel 
 $T$. 
The  quotient $F(\Fq)/I(E(\Fq))$ is isomorphic  to  $T$.
Take    $a$ in  $F(\Fq)$ such that 
$a\bmod I(E(\Fq))$  generates this quotient.
The fiber  $I^{-1}(a)$
is an irreducible divisor. This means that the  $d$ geometric points above
 $a$ are defined on a degree $d$ extension  $\bL$ of
 $\bK$ and permuted by Galois action. We denote by 
$B=I^{-1}(a)$ the corresponding prime divisor.

Since  $\bL$ is the residue extension of $E$ at $B$, we can represent
elements of $\bL$ in the following way: If $f$ is a function on $E$
with polar divisor disjoint to $B$, we denote by $f\bmod B \in \bL$
the residue of  $f$ at $B$.

% Let $X$, $Y$, $Z$ be projective coordinates on $E$.  telles que $Z$ ne
% s'annule pas sur et $x=\frac{X}{Z}$ et $y=\frac{Y}{Z}$ les coordonn{\'e}es
% affines associ{\'e}es. On pose $\alpha = x\bmod B$ et $\beta = y \bmod B$.
% Un {\'e}lement de $\bL$ est donc donn{\'e} comme un polyn{\^o}me en $\alpha$
% et $\beta$.

For  $f$  a function
in $\Fq(E)$, the degree
 of $f$ is the number of poles of $f$ counted
with multiplicities.
For every  $k\ge 0$ we call  $\cF_k$ the set of degree  $\le k$
functions in  $\Fq(E)$, having no pole at  $B$.
We denote by  $V_k$ the corresponding set of residues in   $\bL$,
$$V_k=\{f\bmod B | f \in \cF_k     \}.$$

We have  
$V_0=V_1=\bK\subset V_2\subset \dots \subset V_{d}=\bL$
(Riemann-Roch) and
$V_k\times V_l\subset V_{k+l}.$
It is clear also 
that  $\cF_k$ is Galois invariant since composition by
a translation from $T$ does not affect the degree of a function. 
Therefore $V_k$ is  invariant under the action of  $\Gal(\bL/\bK)$.

If we want to test whether an element
 $z$ of $\bL$ is in  $V_k$, we look for a function
$f$ in  $\cF_k$ such that  $f=z\pmod B$. This is an
interpolation problem which is hardly more
difficult than in the
two previous cases (polynomials for Kummer and
rational fractions for the torus).
We look for $f$ as a quotient
of two homogeneous forms
of degree  $\lceil \frac{k+1}{3}\rceil$, which can be done with linear algebra.

One can choose a smoothness basis
consisting of all elements $f\bmod
B$ in  $V_\kappa$ for a given  $\kappa$.
Factoring an element
 $z=f\bmod B$ of $\bL$ boils down to factoring the divisor
of $f$ as a sum of prime divisors of degree $\le \kappa$.

What conditions are sufficient  for 
an elliptic curve to exist with all the required properties ?
Since the number of $\Fq$-rational points on the elliptic curve
is  divisible by  $d$, the size $q$ of the field cannot be too
small, that is
$$q+2\sqrt q+1>  d.$$

Assume  $d$ is odd
and there exists a squarefree multiple  $D$ of $d$ such that 
 $D\not \equiv 1\bmod p$ and
$$q+1-2\sqrt q < D < q+1+2\sqrt q.$$
There exists an ordinary elliptic curve  $E$  over  $\Fq$ having $D$ 
rational points over $\Fq$ and trace  $t=q+1-D$. 
The ring  $\ZZ[\phi]$  is integrally closed locally at every odd prime dividing
$D$. The  larger ring  $\End(E)$  has the same property.
The ideal  $(\phi -1)$ of  $\End(E)$ has a unique degree $d$ factor $\igot$.
The quotient $\End(E)/\igot$ is cyclic and  $\igot$ is invertible
in  $\End(E)$.

Given $q$ and $\phi$ (a quadratic integer) as above, one can find
an elliptic curve $E/\Fq$ by exhaustive search or using complex multiplication
theory.

\subparagraph*{Example.} Let  $p=q=11$, 
 and 
$d=D=7$, so  $t=5$ and  $\phi^2-5\phi+11=0$.
The elliptic curve  $E$ with equation
$y^2 + xy = x^3 + 2x + 8$
has complex multiplication by $\ZZ[\frac{\sqrt{-19}+1}{2}]$.
The  discriminant of  $\ZZ[\phi]$ is 
$-19$, so  $\End(E)=\ZZ[\phi]$. The ideal 
 $\igot =(\phi-1)$ is invertible and its kernel  $T$ is
the full group of $\Fq$-rational points on $E$.
The kernel of the degree $7$ isogeny  $I: E\rightarrow F$ is the group
of rational points on 
$E$ and
for any non zero   $a\in F(\FF_{11})$, the fiber  $B=I^{-1}(a)$ is  irreducible.

\section{Sieving algorithms and surfaces}\label{section:JL}

There exists a family of algorithms for factoring integers and
computing discrete logarithms that rely on intersection theory on
algebraic or arithmetic surfaces. These algorithms are known
as {\it the number field sieve}, {\it the function field sieve},
 etc. The core of these algorithms is illustrated
on the front page of the book~\cite{nfs}.
In this section, we present the ideas underlying 
this family of algorithms in a rather general setting.
This will help us to describe  our construction  in the next section
\ref{section:carre}. The sieving algorithm 
invented by Joux and Lercier in~\cite{JouxLercier}
for computing  discrete logarithms 
will serve as  a nice illustration for these ideas.

Let  $\Fp$ be the field with  $p$ elements where $p$ is prime.
Let $\cS$ be a smooth projective reduced, absolutely irreducible 
surface over $\Fp$. Let  $\cA$ and  $\cB$ be two absolutely
irreducible curves
on  
$\cS$.
Let  $\cI$ be  an irreducible
sub-variety of the intersection  $\cA\cap\cB$. We assume
that $\cA$ and  $\cB$ meet transversely at $\cI$ and we
denote by  $d$ the degree of  $\cI$.
The residue field of  $\cI$ is  $\Fp(\cI)=\Fq$ with $q=p^d$.

We need a pencil (linear or at least algebraic and connected)
of effective divisors on $\cS$. We denote it by
 $(D_\lambda)_{\lambda \in \Lambda}$ 
where $\Lambda$ is the parameter space.

We fix an integer  $\kappa$ and we look (at random) for divisors $D_\lambda$ 
in the pencil, such that both intersection divisors  $D\cap \cA$ 
and  $D\cap \cB$
are  disjoint to  $\cI$ and  $\kappa$-smooth (they split
as sums of effective $\Fq$-divisors of degree $\le \kappa$).

We define an  equivalence 
relation
 $\equiv_\cI$  on the set of  divisors on 
$\cS$ not meeting  $\cI$: 
We say  $D\equiv_\cI 0$  if and only if  $D$ is the divisor of a function
 $f$ and $f$ is constant  modulo $\cI$.
The equivalence classes for this relation are parameterized by points 
in  some  algebraic group denoted  $\Pic(\cS,\cI)$. This algebraic group is
an extension
of  $\Pic(\cS)$ by a torus $T_\cI$ of  dimension $d-1$.

One similarly defines the algebraic groups  $\Pic(\cA,\cI)$ and 
$\Pic(\cB,\cI)$. These are generalized jacobians of  $\cA$ and  $\cB$ respectively.
The natural (restriction)  morphisms $\Pic(\cS,\cI)\rightarrow \Pic(\cA,\cI)$
and 
$\Pic(\cS,\cI)\rightarrow \Pic(\cB,\cI)$ induce the identity on the torus $T_\cI$.

Let  $N$ be an integer that kills the three groups 
$\Pic^0(\cS)(\Fp)$, $\Pic^0(\cA)(\Fp)$, and  $\Pic^0(\cB)(\Fp)$.
Let  $\lambda$ and  $\mu$ be two parameters in  $\Lambda$
corresponding to the divisors $D_\lambda$ and $D_\mu$ in our pencil. We assume
that  $D_\lambda\cap \cA$, $D_\mu \cap \cA$,
$D_\lambda\cap \cB$, and  $D_\mu \cap \cB$ are
smooth and disjoint to  $\cI$.

Let  $D_\lambda\cap \cA = \sum {\mathfrak a}_i$,  $D_\mu\cap \cA = \sum {\mathfrak b}_j$, 
$D_\lambda\cap \cB = \sum {\mathfrak c}_k$ and  $D_\mu\cap \cB = \sum
{\mathfrak d}_l$ be decompositions as sums of effective
divisors on  $\cA$ and $\cB$ with degree  $\le \kappa $.
The divisor  $D_\lambda-D_\mu$ is algebraically equivalent to zero
and $N(D_\lambda-D_\mu)$ is principal.

Let  $f$ be a  function on  $\cS$ with divisor  $N(D_\lambda-D_\mu)$.
We fix a smooth divisor  $X$ on  $\cA$ (resp. $Y$ on $\cB$) with degree  $1$.
For every  $i$ and $j$, let  $\alpha_i$ and $\beta_j$ be  functions on  $\cA$ 
with divisors  $N({\mathfrak a}_i-\deg({\mathfrak a}_i)X)$ and  $N({\mathfrak b}_j-\deg({\mathfrak b}_j)X)$.
Similarly, for every  $k$ and $l$, let   $\gamma_k$ and $\delta_l$  be functions on  $\cB$ 
with divisors  $N({\mathfrak c}_k-\deg({\mathfrak c}_k)Y)$ and $N({\mathfrak d}_l-\deg({\mathfrak d}_l)Y)$.
There exist two  multiplicative 
constant $c$ and $c'$ in $\Fps$ such that  
$$f\equiv c . \frac{\prod_i\alpha_i}{\prod_j\beta_j}  \equiv  c' . \frac{\prod_k\gamma_k}{\prod_l\delta_l}\bmod \cI.$$

This congruence 
 can be regarded as  a relation in the group  $T_\cI(\Fp)=\Fqs/\Fps$.
The factors in the first fraction  belong to the smoothness basis on the $\cA$
side: They are residues modulo $\cI$ of functions on $\cA$ with degree $\le
\kappa$. Similarly, 
the factors in the second fraction  belong to the smoothness basis on the $\cB$
side: They are residue modulo $\cI$ of functions on $\cB$ with degree $\le
\kappa$.

Joux and Lercier take $\cS/\Fp$ to be  $\cS=\PU \times \PU$ the
product of $\PU$ with itself over $\Fp$. To avoid
any confusion we call  $\cC_1=\PU/\Fp$ the first factor and  
$\cC_2=\PU/\Fp$ the second factor. 
Let  $O_1$ be a rational  point on $\cC_1$ and  $\cU_1=\cC_1-O_1$. Let 
$x$ be an affine coordinate on  $\cU_1\sim \AU$. We similarly choose   $O_2$, $\cU_2$
and $y$ an affine coordinate on $\cU_2$.

They choose  $\cA$ to be the Zariski closure in  $\cS$ of the
curve in
 $\cU_1\times \cU_2$
with equation  $y=f(x)$ where  $f$ is a polynomial with degree $d_f$ in $\Fp[x]$.
As for  $\cB$, they choose the Zariski closure in  $\cS$ of the curve with
equation  $x=g(y)$ where  $g$ is a  polynomial with degree $d_g$ in  $\Fp[y]$.

The N{\'e}ron-Severi group of a product of two smooth
algebraically irreducible
 projective
curves is $\ZZ$ times $\ZZ$ times the group of homomorphisms
between the jacobians of the two curves. See
\cite[Mumford's appendix to Chapter VI]{zar}. The Hurwitz  formula for the intersection of two
classes is also given in this appendix.

Here the
   N{\'e}ron-Severi  group  of  $\cS$ is $\ZZ\times \ZZ$. The algebraic
equivalence  class of a divisor  $D$ is given as its bidegree 
$(\dx(D),\dy(D))$
where  $\dx(D)=D.(\cC_1\times O_2)$ and   $\dy(D)=D.( O_1\times
\cC_2)$. The intersection form is given by the formula 
$$D.E=\dx(E)\dy(D)+\dx(D)\dy(E).$$

The bidegree of $\cA$ is  $(d_f,1)$ and the bidegree of  $\cB$ is $(1,d_g)$.
So  $\cA.\cB=1+d_fd_g$ and the intersection of  $\cA$ and  $\cB$ is made of
the point $O_1\times O_2$ and the  $d_fd_g$ points of the  form
$(\alpha,f(\alpha))$ where 
$\alpha$
is one of the  $d_fd_g$  roots of  $g(f(x))-x$.

Let $h(x)$ be a simple  irreducible factor of the later polynomial and let $d$
be its degree.
We take  $\cI$ to be the zero dimensional and degree $d$ corresponding variety.
The residue field  $\Fp(\cI)$ is finite of order $q$ where  $q=p^d$.

To finish with, we need a pencil of effective divisors  $(D_\lambda)_{\lambda \in \Lambda}$ 
on  $\cS$. It is standard to take for $\Lambda$ the set of polynomials
$\lambda$ in 
$\Fp[x,y]$  with given bidegree   $(u_x,u_y)$ where $u_x$ and $u_y$ are chosen
according to $p$ and $q$.  The corresponding divisor  $D_\lambda$
to  $\lambda$ 
is the Zariski  closure of the zero set of  $\lambda$. It has bidegree  
 $(u_x,u_y)$ too.

We fix an integer  $\kappa$
and look for divisors $D_\lambda$ such that
the two  intersection divisors  $D_\lambda\cap \cA$ and  $D_\lambda \cap \cB$
are  disjoint to  $\cI$ and  $\kappa$-smooth.
For example, if  $\lambda (x,y)$ is a polynomial in $x$ and  $y$, the
intersection of  $D_\lambda$ and  $\cA$ has degree  $d_fu_y+u_x$. Its affine
part is given by the roots of the polynomial
$\lambda(x,f(x))=0$.
Similarly, the intersection of  $D_\lambda$ and  $\cB$ has degree  $u_y+u_xd_g$. Its
affine part is given by the roots of the polynomial 
$\lambda(g(y),y))=0$.
Joux and Lercier explain how to choose  $u_x$,  $u_y$ and $\kappa$ according
to  $p$ and  $d$.

\section{Finite residue fields on elliptic squares}\label{section:carre}

In this section we try to conciliate
the generic construction in 
section~\ref{section:JL} and the ideas developed in 
section~\ref{section:ell}.
We would like the   automorphisms of  $\Fp(\cI)$ to be induced by 
automorphisms
of the surface  $\cS$.
So let 
$E$ be an ordinary  elliptic curve over  $\Fp$ 
and let  $\igot$ be an  invertible  ideal in the endomorphism ring
$\End(E)$. 
We assume  that $\igot$ divides   $\phi-1$ and  $\End(E)/\igot$ is cyclic
of order $d\ge 2$.
So   $E(\Fq)$ contains 
a cyclic subgroup  $T=\Ker \igot$ of order  $d$.
Let $I: E \rightarrow F$ be the quotient by   
$\Ker \igot$ isogeny
and let  $J: F\rightarrow E$ be such that  $\phi-1=J\circ I$.

We take for  $\cS$ the product  $E\times E$ and to avoid any confusion,
we call $E_1$ the first  factor and  $E_2$ the second  factor. Let  $O_1$ be
the origin on  $E_1$ and  $O_2$ the origin on  $E_2$. 

We use
again the description
of the  N\'eron-Severi group of a product of two curves as
given in \cite[Appendix to Chapter VI]{zar}. 
This time, the N{\'e}ron-Severi group of $\cS$ is  $\ZZ\times \ZZ \times \End(E)$. The
class  $(d_1,d_2,\xi)$ of a divisor $D$ consists of the bidegree
and the induced isogeny. More precisely,  $d_1$
is the intersection degree of  $D$ and  $E_1\times O_2$, $d_2$ is the
intersection degree of 
$D$ and $O_1\times E_2$, and  $\xi$ is the homomorphism from  $E_1$ to $E_2$
induced
by the correspondence associated with $D$.

Let  $\alpha$ and  $\beta$ be two endomorphisms of  $E$ and let 
$a$ and  $b$ be two   $\Fp$-rational  points on  $E$. 
We take   $\cA$ to be the inverse image of  $a$ by the morphism from  $E\times E$
to $E$ that maps  $(P,Q)$ onto  $\alpha(P)-Q$.
Let  $\cB$ be the inverse image 
of  $b$ by the morphism from $E\times E$
onto  $E$ that sends  $(P,Q)$ onto  $P-\beta(Q)$.

Assume  $1-\beta\alpha=\phi-1$.  The intersection 
of  $\cA$ and  $\cB$
consists of points  $(P,Q)$ such that $(\phi-1)(P)=b-\beta(a)$ and 
$Q=\alpha(P)-a$.

We choose  $a$ and  $b$ such that there exists a point
 $c$ in  $F(\Fp)$ generating  $F(\Fp)/I(E(\Fp))$ and
satisfying $J(c)=b-\beta(a)$. Then the intersection
between  $\cA$ and  $\cB$ contains an irreducible  component   $\cI$
of degree $d$.

The class of  $\cA$ is  $(\alpha\bar\alpha, 1,  \alpha)$.
Indeed, the first coordinate  of this triple is the degree of 
the projection $\cA\rightarrow E_2$
onto  the second component, that is the number of
solutions in $P$ to 
$\alpha(P)=Q+a$ for generic $Q$. This is the degree $\alpha\bar\alpha$
of $\alpha$.
The second  coordinate  of this triple is the degree of the projection
$\cA\rightarrow E_1$
onto  the first   component, that is the number of
solutions in $Q$ to 
$Q=\alpha(P)-a$ for generic $P$. This is $1$.
The third coordinate is the morphism in $\Hom(E_1,E_2)$ induced 
by the correspondence $\cA$. This is clearly $\alpha$.
In  the same way, we prove that
the  class of $\cB$ is  $(1,\beta\bar\beta, \bar \beta)$.

Now let  $D$ be a divisor on  $\cS$ and  $(d_1,d_2,\xi)$
its class in the  N{\'e}ron-Severi group.
The intersection degree of  $D$ and  $\cA$ is thus 
\begin{equation}\label{eq:int1}
D.\cA=d_1+d_2\alpha\bar\alpha  -\xi\bar \alpha-\bar \xi\alpha
\end{equation}
and similarly 
\begin{equation}\label{eq:int2}
D.\cB=d_1\beta\bar\beta+d_2  -\xi\bar \beta-\bar \xi\beta.
\end{equation}

We are particularly interested in the case where
$\alpha$ and  $\beta$  have norms of essentially the same size 
(that is the square root of the norm of $\phi -2$).
We then obtain a similar behavior as
the algorithm in section 
\ref{section:JL} with an extra advantage: The smoothness
bases on both 
$\cA$ and $\cB$ are  Galois invariant.

Indeed, let  $f_\cA$ be a function with degree  $\le \kappa$ on 
$\cA$. A point on  $\cA$ is  a couple $(P,Q)$ with 
$Q=\alpha(P)-a$. So the projection on the first component 
$\Pi_1 : E_1\times E_2\rightarrow E_1$ is an isomorphism. 
There is a unique function $f_1$ on $E_1$ such that $f_\cA = f_1\circ \Pi_1$.
Assume now 
that  $(P,Q)$ is in  $\cI\subset \cA$. Then $f_\alpha (P,Q)=f_1(P)$ is an element
of the smoothness basis on $\cA$. We observe that 
 $f_1(P)^p=f_1(\phi(P))=f_1(P+t)$ where  $t$ is in the kernel  $T$ of $\igot$. 
So $f_1(P)^p$ is the value at  $P$ of 
$f_1\circ\tau_{t}$ where  $\tau_t: E_1\rightarrow E_1$
is the  translation by $t$.
Since  $f_1\circ \tau_t$ and  $f_1$ have the same degree,
the value of $f_1\circ \tau_t$ at  $P$ is again an element
in the smoothness basis.

That way, one can divide by  $d$ the size of
either smoothness basis on  $\cA$ and $\cB$.

As in section \ref{section:JL} we need a pencil  of divisors 
on $\cS$ with small class in the N{\'e}ron-Severi group. 
We choose small values for  $(d_1,d_2,\xi)$ that minimize
the expressions in 
 Eq.~(\ref{eq:int1}) and  Eq.~(\ref{eq:int2})
under the three contraints 
    $d_1\ge 1$, $d_2\ge 1$ and 
\begin{equation}\label{eq:ine}
d_1d_2\ge \xi\bar\xi+1. 
\end{equation}

We look for effective divisors in the algebraic equivalence   class $\cgot = (d_1,d_2,\xi)$. 
Recall   $O_1$ is the origin on  $E_1$  and  $O_2$ the origin on  $E_2$.
The  graph $\cG = \{(P,Q) | Q=-\xi (P)  \}$ of  $-\xi: E_1\rightarrow E_2$ 
is a divisor in the class $(\xi\bar\xi , 1, -\xi)$.
The divisor  
$\cH=-\cG+(d_1+\xi\bar\xi )O_1\times E_2+(d_2+1)E_1\times O_2$ 
is in  $\cgot$.
We compute the linear space 

$$\cL(-\cG+(d_1+\xi\bar\xi )O_1\times
E_2+(d_2+1)E_1\times O_2)$$
\noindent  using the (restriction) exact sequence
\begin{multline*}
0\rightarrow \cL_S(-\cG+(d_1+\xi\bar\xi)O_1\times E_2
+(d_2+1)E_1\times O_2)\\
\rightarrow \cL_{E_1}((d_1+\xi\bar\xi )O_1)\otimes \cL_{E_2}((d_2+1) O_2)
\rightarrow \cL_\cG(\Delta)
\end{multline*}
where  $\Delta$ is the divisor on  $\cG$ given by the intersection
with 
$$(d_1+\xi\bar\xi)O_1\times E_2+(d_2+1)E_1\times O_2.$$

This divisor has degree
 $d_1+\xi\bar\xi +(d_2+1)\xi\bar\xi$, so the 
dimension
of the right hand term in the sequence above
is equal to this number.

On the other hand, the middle term has dimension $(d_1+\xi\bar\xi)(d_2+1)$,
that is strictly bigger than the dimension of the right hand term, because of
Inequality~(\ref{eq:ine}).
So the linear space on the left is
non zero and the divisor class is effective. Inequality~(\ref{eq:ine}) is a sufficient condition for
effectivity.

In practice, one computes a basis for
$\cL_{E_1}((d_1+\xi\bar\xi)O_1)$
and a basis for  $\cL_{E_2}((d_2+1) O_2)$ and one multiplies
the two basis  (one takes all products
of one element in the first basis with one element in the second
basis.) This produces
a basis for $\cL_{E_1}((d_1+\xi\bar\xi )O_1)\otimes \cL_{E_2}((d_2+1) O_2)$.

One selects enough   (more than  
$d_1+\xi\bar\xi+(d_2+1)\xi\bar\xi$)  points 
$(A_i)_i$ on 
$\cG$ and  one evaluates all functions
in the above basis at all these points.
A linear algebra calculation
produces  a basis for the subspace
of $\cL_{E_1}((d_1+\xi\bar\xi )O_1)\otimes \cL_{E_2}((d_2+1) O_2)$
consisting of functions that vanish along $\cG$.
For every function $\phi$ in the later subspace, the divisor 
of zeroes of  $\phi$ contains  $\cG$ and the difference
$(\phi)_0
-\cG$ is an effective divisor in the
linear equivalence class of $\cH$.

We have thus constructed a complete
linear equivalence class inside $\cgot$.
To find the other linear classes
in  $\cgot$, we remind that   $E\times E$ is isomorphic
to its Picard variety. So it suffices  to replace $\cH$ 
in the above calculation by 
$\cH+E_1\times Z_2-E_1\times O_2+Z_1\times E_2-O_1\times E_2$ 
where  $Z_1$ and  $Z_2$ run over  $E_1(\Fp)$ and  $E_2(\Fp)$ respectively.

\section{Experiments}\label{section:exp}

In this section, we give a practical example of the geometric construction
of section~\ref{section:carre}.
We 
perform a discrete logarithm computation in $\FF_{61^{19}}$. 
In such a field,
Joux and Lercier algorithm would handle a factor basis of irreducible
polynomials of degree 2 over $\FF_{61}$, in two variables. Such a factor
basis would have  about $3600$ elements. It turns out that in this case we can reduce
the factor basis to only $198$ elements using the ideas given in  the previous section.

\paragraph*{Initialization phase.}

We set $p=61$ and consider the plane projective elliptic curve $E$ over $\Fp$
with equation $Y^2Z = X^3 + 20XZ^2 + 21Z^3$.  It is ordinary with trace
$t=-14$. The ring generated by the Frobenius $\phi$ has discriminant
$-48$. The full endomorphism ring of $E$ 
is the maximal order in the field $\QQ(\sqrt{-3})$.

Let  $\beta$ be the degree $3$ endomorphism of $E$ given by 
\begin{displaymath}
  \begin{array}{crcl}
    \beta: &E & \rightarrow  & E\,, \\
&    (x~:y~:1) & \mapsto & (\frac
{20\,{x}^{3}+36\,{x}^{2}+35\,x+40}{\left( x+7 \right) ^{2}}~:y\frac {
  58\,{x}^{3}+59\,{x}^{2}+12\,x+21}{\left( x+7 \right) ^{3}}~: 1)\,.
  \end{array}
\end{displaymath}

We check $\beta^2=-3$ and we fix an isomorphism between $\End(E)\otimes \QQ$
and $\QQ(\sqrt{-3})\subset \CC$  by setting $\beta = \sqrt{-3}$. The Frobenius
endomorphism is $\phi=-7+2\sqrt{-3}$.

Let  $\alpha$  be the degree $4$  endomorphism  
defined by $\alpha= 1+\beta=1+\sqrt{-3}$. It can be  given explicitly by 
\begin{displaymath}
  \begin{array}{crcl}
\alpha:&    E & \rightarrow  & E\,, \\
&    (x~:y~:1) & \mapsto & ({\frac
        {49\,{x}^{4}+28\,{x}^{3}+55\,{x}^{2}+53\,x+27}{ \left( x+25 \right)
          \left( x+27 \right) ^{2}}}~:y\frac {
          38\,{x}^{5}+37\,{x}^{4}+30\,{x}^{3}+49\,{x}^{2}+9\,x+46
       }{\left( x+25 \right) ^{2} \left( x+27 \right) ^{3}}~: 1 )\,.
  \end{array}
\end{displaymath}

The endomorphism  $I=1-\beta\alpha$ has  degree 19 and divides $\phi-1$.
The kernel of $I$ consists of the following $19$ rational points,
\begin{footnotesize}
\begin{multline*}
  \Ker I = \{ (0: 1: 0), (11: \pm 13: 1), (14: \pm19: 1), (21: \pm8:
  1), (35: \pm15: 1), \\
  (40: \pm10: 1), (41: \pm10: 1), (45: \pm27: 1), (48: \pm2: 1), (51: \pm23: 1) \}\,.
\end{multline*}  
\end{footnotesize}
Let $\cS=E\times E$.  We  call $E_1=E$ the first factor
and $E_2=E$ the second one.
If  $P$ and $Q$ are independent generic points on $E$,
then $(P,Q)$ is a generic point on $\cS$.
Let $a$ on $E$
be the point with coordinates $(52:24:1)$. 
Let $\cA\subset \cS$ be the curve with equation
$\alpha(P) - Q = a$. 
Let $b$ on $E$
be the point with coordinates $(1:46:1)$. 
Let $\cB\subset \cS$ be the curve with equation
$P - \beta(Q) = b$.  The numerical class of $\cA$
is $(4,1,1+\sqrt{-3})$ and the numerical class of
$\cB$ is $(1,3,-\sqrt{-3})$.
Note that $b-\beta(a)=(57:11:1)$ is of order $38$ and generates $E(\Fp)$
modulo the image of $I$.

Call $\cI$ the intersection $\cA\cap \cB$. It consists 
of points  $(P,Q)$ such that $(1-\beta\alpha)(P)=b-\beta(a)$,
$Q=\alpha(P)-a$ and thus $(\alpha\beta-1)(Q) = a - \alpha(b)$.
In terms of the affine coordinates $(x_1,y_1)$ of $P$ and $(x_2,y_2)$ of $Q$,
this reads 
\begin{footnotesize}
\begin{multline}\label{eqxp}
  x_1 = {\frac { \left( 44\,{x_2}^{4}+12\,{x_2}^{3}+9\,{x_2}^{2}+46\,x_2+40 \right) y_2
}{{x_2}^{6}+34\,{x_2}^{5}+41\,{x_2}^{4}+47\,{x_2}^{3}+7\,{x_2}^{2}+14\,x_2+58}}+\\{
\frac {{x_2}^{6}+26\,{x_2}^{5}+25\,{x_2}^{3}+41\,{x_2}^{2}+19\,x_2+6}{{x_2}^{6}+34
\,{x_2}^{5}+41\,{x_2}^{4}+47\,{x_2}^{3}+7\,{x_2}^{2}+14\,x_2+58}},
\end{multline}
\begin{multline}\label{eqyp}
  y_1 = {\frac {
 \left( 11\,{x_2}^{7}+2\,{x_2}^{6}+50\,{x_2}^{5}+59\,{x_2}^{4}+57\,{x_2}^{3}+30
\,{x_2}^{2}+4\,x_2+14 \right) y_2}{{x_2}^{9}+51\,{x_2}^{8}+7\,{x_2}^{7}+32\,{x_2}^{6
}+56\,{x_2}^{5}+48\,{x_2}^{4}+26\,{x_2}^{3}+49\,{x_2}^{2}+18\,x_2+41}}+\\{\frac {
46\,{x_2}^{9}+54\,{x_2}^{8}+2\,{x_2}^{7}+4\,{x_2}^{6}+52\,{x_2}^{5}+17\,{x_2}^{4}+
60\,{x_2}^{3}+41\,{x_2}^{2}+48\,x_2+21}{{x_2}^{9}+51\,{x_2}^{8}+7\,{x_2}^{7}+32\,{
x_2}^{6}+56\,{x_2}^{5}+48\,{x_2}^{4}+26\,{x_2}^{3}+49\,{x_2}^{2}+18\,x_2+41}}\,,
\end{multline}
\end{footnotesize}
or alternatively, $x_2$, $y_2$ can be given as functions of degree $8$ and degree $12$ in
$x_1$, $y_1$.

The projection of $\cI$ on $E_1$ (resp. $E_2$) yields a place $\PP$ (resp.
$\PQ$) of degree 19 defined in the affine coordinates $(x,y)$ by the equations
\begin{footnotesize}
\begin{multline*}
  \PP = ({x_1}^{19}+60\,{x_1}^{18}+25\,{x_1}^{17}+21\,{x_1}^{16}+23\,{x_1}^{15}+22\,{x_1}^
{14}+49\,{x_1}^{13}+38\,{x_1}^{12}+30\,{x_1}^{11}+57\,{x_1}^{10}+\\
3\,{x_1}^{9}+15
\,{x_1}^{8}+26\,{x_1}^{7}+17\,{x_1}^{6}+45\,{x_1}^{5}+30\,{x_1}^{4}+48\,{x_1}^{3}+
55\,{x_1}^{2}+18\,x_1+35,\\
y_1+12\,{x_1}^{18}+38\,{x_1}^{17}+5\,{x_1}^{16}+{x_1}^{15}+
45\,{x_1}^{14}+42\,{x_1}^{13}+18\,{x_1}^{12}+34\,{x_1}^{11}+39\,{x_1}^{10}+\\59\,{
x_1}^{9}+16\,{x_1}^{8}+18\,{x_1}^{7}+16\,{x_1}^{6}+36\,{x_1}^{5}+11\,{x_1}^{4}+9\,
{x_1}^{3}+48\,{x_1}^{2}+59\,x_1+8)\,,
\end{multline*}
\begin{multline*}
  \PQ = ({x_2}^{19}+25\,{x_2}^{18}+34\,{x_2}^{17}+46\,{x_2}^{16}+16\,{x_2}^{15}+14\,{x_2}^
{14}+58\,{x_2}^{13}+52\,{x_2}^{12}+39\,{x_2}^{11}+48\,{x_2}^{10}+\\18\,{x_2}^{9}+
56\,{x_2}^{8}+41\,{x_2}^{7}+40\,{x_2}^{6}+11\,{x_2}^{5}+33\,{x_2}^{4}+55\,{x_2}^{3
}+14\,{x_2}^{2}+5\,x_2+56,\\y_2+42\,{x_2}^{18}+40\,{x_2}^{17}+23\,{x_2}^{16}+41\,{x_2}
^{15}+14\,{x_2}^{14}+12\,{x_2}^{13}+30\,{x_2}^{12}+50\,{x_2}^{11}+33\,{x_2}^{10}
+\\33\,{x_2}^{9}+60\,{x_2}^{8}+15\,{x_2}^{7}+54\,{x_2}^{6}+13\,{x_2}^{5}+17\,{x_2}^{
4}+31\,{x_2}^{3}+50\,{x_2}^{2}+52\,x_2+3)\,.
\end{multline*}
\end{footnotesize}
The residue fields of these two places are isomorphic (both being degree
$19$ extensions of $\FF_{61}$).  
We fix an  isomorphism between these two residue fields by setting 
\begin{footnotesize}
\begin{multline}\label{eqiso}
  x_2 \mapsto 2\,{x_1}^{18}+57\,{x_1}^{17}+21\,{x_1}^{16}+10\,{x_1}^{15}+54\,{x_1}^{14}+35\,{x_1
  }^{13}+45\,{x_1}^{12}+27\,{x_1}^{11}+41\,{x_1}^{10}+\\
55\,{x_1}^{9}+27\,{x_1}^{8}+
  36\,{x_1}^{7}+29\,{x_1}^{6}+50\,{x_1}^{5}+44\,{x_1}^{4}+18\,{x_1}^{3}+38\,{x_1}^{2
  }+51\,x_1+18\,.
\end{multline}
\end{footnotesize}
Fixing this isomorphism is equivalent to choosing a geometric point in $\cI$.

\paragraph*{Sieving phase.}

We are now going to look for ``smooth'' functions on $\cS$. We first explain
what  we mean by smooth in this context.
Let $\varepsilon(x_1, y_1, x_2, y_2)$ be a function on $\cS$. We assume
$\varepsilon$ does not vanish at $\cI$.
Let $\Pi_1: \cS=E_1\times E_2\rightarrow E_1$ be the  projection on the first
factor. The restriction of $\Pi_1$ to $\cA$ is a bijection. So we can 
define a point on $\cA$ by its coordinates $(x_1,y_1)$.
Let $\Pi_2: \cS=E_1\times E_2\rightarrow E_2$ be the  projection on the second
factor. The restriction of $\Pi_2$ to $\cB$ is a bijection. So we can 
define a point on $\cB$ by its coordinates $(x_2,y_2)$.

Let $\vare_1(x_1, y_1)$ (resp.  $\varepsilon_2(x_2, y_2)$) be the restriction
of $\vare$ to $\cA$ (resp. $\cB$). For example $\vare_2(x_2,y_2)$ is obtained
by substituting $x_1,y_1$ as functions in $x_2,y_2$ in $\varepsilon$ thanks to
Eq.~(\ref{eqxp}) and Eq.~(\ref{eqyp}).

 The function 
$\varepsilon$ is said to be smooth if the divisors of $\varepsilon_1$ and
$\varepsilon_2$ both contain only places of small degree $\kappa$. In our example, we
choose $\kappa=2$. Let us remark at this point that thanks to the isomorphism given
by Eq.~(\ref{eqiso}), the reduction modulo $\PP$ of $\varepsilon_1$ is equal
to the reduction modulo $\PQ$ of $\varepsilon_2$, and this yields an equality
in $\FF_{61^{19}}$.

To every non-zero function on $\cS$, one can associate a linear pencil of
divisors. We define  the linear (resp. numerical) class of the function
to be the linear (resp. numerical) class of the divisor
of its zeroes (or poles).

We shall be firstly interested in functions $\varepsilon$ with numerical class
$(1, 0, 0)$. An effective divisor in these classes is $c\times
E_2$ where $c$ is a place of degree $1$ on $E_1$ and it is not
difficult to see that the intersection degrees of such a divisor with
$\mathcal A$ and $\mathcal B$ are $1$ and $3$. Functions with numerical class
$(2, 0, 0)$ are obtained in the same way.

We found similarly functions $\varepsilon$ in the class $(0,1,0)$,
derived from divisors $E_1\times c$. The intersection degrees are now
$4$ and $1$. Functions with numerical class $(0, 2, 0)$ are obtained
in the same way too.  More interesting, the class $(1, 1, 1)$
 containing the divisors with equation $P=Q+c$, yields
intersection degrees $3$ and $4$.

We finally consider the class $(2,2,1)$ which is, by far, much larger than the
previous classes. The intersection degrees are $8$ and $8$. To
enumerate functions in this class, we first build a basis for  the linear space
associated to  divisors of degree $3$ on both  $E_1$ and $E_2$. For instance, let us consider
${\mathcal L}_{E_1}(3\,O_1)$ and ${\mathcal L}_{E_2}(3\,O_2)$,  basis of
which are given by $\{1, x_1, y_1\}$ and $\{1, x_2, y_2\}$. We then
determinate that a basis for the subspace of ${\mathcal L}_{E_1}(3\,O_1)\otimes
{\mathcal L}_{E_2}(3\,O_2)$, consisting of functions that vanish along the
graph ${\mathcal G} = \{ (P,Q), Q = -P\}$, is given by $\{y_1\,x_2+x_1\,y_2,
y_1+y_2, x_1-x_2\}$. An exhaustive enumeration of functions of the form
$y_1\,x_2+x_1\,y_2 +\lambda( y_1+y_2) + \mu ( x_1-x_2)$, with $\lambda,\mu \in
\FF_{p}$ yields useful equations.

We give examples of such relations in Tab.~\ref{tab:relex}.

\begin{table}[htbp]
  \centering
  \begin{footnotesize}
  \begin{tabular}{cp{0.42\textwidth}p{0.42\textwidth}}
    Class & \multicolumn{1}{c}{div $\varepsilon_1$} & \multicolumn{1}{c}{div
      $\varepsilon_2$} \\\hline
    $(1,0,0)$ &
    $(x_1 + 43, y_1 + 33) - (x_1 + 13, y_1 + 59)$ &
    $(x_2^2 + x_2 + 52, y_2 + 10\,x_2 + 37) + (x_2 + 
    12, y_2 + 35) - (x_2 + 2, y_2 + 20) - (x_2^2 + 26\,x_2 + 39, y_2 + 5\,x_2
    + 27)$\\\hline
    $(2,0,0)$ &
    $(x_1^2 + 56\,x_1 + 34, y_1 + 22\,x_1 + 52) - 2\,(x_1 + 13, y_1 + 59)$ &
    $(x_2^2 + 37\,x_2 + 53, y_2
    + 42\,x_2 + 58) + (x_2^2 + 12\,x_2 + 19, y_2 + 52\,x_2 + 43) + (x_2^2 + 41\,x_2 + 29, y_2 + 
    33\,x_2 + 41) - 2\,(x_2 + 2, y_2 + 20) - 2\,(x_2^2 + 26\,x_2 + 39, y_2 + 5\,x_2 + 27)$\\\hline
    $(0,1,0)$ & $
    (x_1^2 + 4\,x_1 + 12, y_1 + 55\,x_1 + 47) + (x_1^2 + 45\,x_1 + 31, y_1 + 19\,x_1 + 23) - (x_1 + 
    42, y_1 + 60) - (x_1 + 36, y_1 + 15) - (x_1^2 + 60\,x_1 + 25, y_1 +
    36\,x_1 + 26)$ & $(x_2 + 
    43, y_2 + 33) - (x_2 + 13, y_2 + 59)$ \\\hline
    $(0,2,0)$ & $
    (x_1^2 + 26\,x_1 + 12, y_1 + 12\,x_1 + 32) + (x_1^2 + 48\,x_1 + 6, y_1 + 59) + (x_1^2 + 53\,x_1 
    + 56, y_1 + 42\,x_1 + 56) + (x_1^2 + 3\,x_1 + 38, y_1 + 17\,x_1 + 36) - 2\,(x_1 + 42, y_1 + 
    60) - 2\,(x_1 + 36, y_1 + 15) - 2\,(x_1^2 + 60\,x_1 + 25, y_1 +
    36\,x_1 + 26)$ & $ (x_2^2 + 
    24\,x_2 + 39, y_2 + 37\,x_2 + 27) - 2\,(x_2 + 13, y_2 + 59)$ \\\hline
    $(1,1,1)$ & $
    (x_1 + 2, y_1 + 41) + (x_1^2 + 26\,x_1 + 39, y_1 + 56\,x_1 + 34) - (x_1^2 + 48\,x_1 + 6, y_1 + 2) -
    (x_1 + 52, y_1 + 25) $ & $ (x_2 + 17, y_2 + 21) + (x_2^2 + 57\,x_2 + 11, y_2 + 33\,x_2) + (x_2 + 55, y_2 + 
    33) - (x_2^2 + 49\,x_2 + 42, y_2 + 26) - (x_2^2 + 3\,x_2 + 4, y_2 +
    30\,x_2 + 20)$ \\\hline
    $(2,2,2)$ & $
    (x_1^2 + 25\,x_1 + 42, y_1 + 5\,x_1 + 13) + (x_1^2 + 30\,x_1 + 19, y_1 + 52\,x_1 + 42) + (x_1^2 + 
    59\,x_1 + 30, y_1 + 28\,x_1 + 22) - 2\,(x_1^2 + 48\,x_1 + 6, y_1 + 2) - 2\,(x_1 + 52, y_1 + 25) $ & $ 
    (x_2^2 + 30\,x_2 + 21, y_2 + 50\,x_2 + 52) + (x_2^2 + 41\,x_2 + 8, y_2 + 54\,x_2 + 58) + (x_2^2 + 32\,x_2
    + 20, y_2 + 34\,x_2 + 28) + (x_2^2 + 42\,x_2 + 49, y_2 + 29\,x_2 + 51) - 2\,(x_2^2 + 49\,x_2 + 
    42, y_2 + 26) - 2\,(x_2^2 + 3\,x_2 + 4, y_2 + 30\,x_2 + 20)$ \\\hline
    $(2,2,1)$ & $
    (x_1 + 24, y_1 + 33) + (x_1 + 25, y_1) + (x_1 + 35, y_1) + (x_1 + 60, y_1 + 46) + 
    (x_1^2 + 33\,x_1 + 43, y_1 + 3\,x_1 + 34) + (x_1^2 + 53\,x_1 + 53, y_1 + 24\,x_1 + 33) -
    (x_1 + 1, y_1) - (x_1 + 54, y_1 + 4) - (x_1^2 + 17\,x_1 + 19, y_1 + 41\,x_1 + 21) - 
    (x_1^2 + 51\,x_1 + 53, y_1 + 44\,x_1 + 31) - (x_1^2 + 55\,x_1 + 38, y_1 + 38\,x_1 + 58) $ & $
    (x_2 + 3, y_2 + 42) + (x_2^2 + 7\,x_2 + 20, y_2 + 33\,x_2 + 46) + (x_2^2 + 38\,x_2 + 
    12, y_2 + 58\,x_2 + 6) + (x_2^2 + 42\,x_2 + 35, y_2 + 7\,x_2 + 41) - (x_2 + 1, y_2) -
    (x_2 + 11, y_2 + 42) - (x_2 + 16, y_2 + 34) - (x_2^2 + 26\,x_2 + 12, y_2 + 49\,x_2 + 29)
    - (x_2^2 + 47\,x_2 + 5, y_2 + 7\,x_2 + 14)$ \\\hline
    $(2,2,1)$ & $
    (x_1 + 10, y_1 + 23) + (x_1 + 20, y_1 + x_1 + 30) + (x_1 + 29, y_1 + 1) + (x_1 + 41,
    y_1 + x_1 + 33) + (x_1^2 + 6\,x_1 + 17, y_1 + 25\,x_1 + 16) + (x_1^2 + 25\,x_1 + 12, y_1 
    + 25\,x_1 + 47) - (x_1 + 1, y_1) - (x_1 + 54, y_1 + 4) - (x_1^2 + 17\,x_1 + 19, y_1 
    + 41\,x_1 + 21) - (x_1^2 + 51\,x_1 + 53, y_1 + 44\,x_1 + 31) - (x_1^2 + 55\,x_1 + 
    38, y_1 + 38\,x_1 + 58)$ & $ (x_2 + 29, y_2 + 60) + (x_2 + 36, y_2 + 15) + (x_2^2 + 
    15\,x_2 + 58, y_2 + 41\,x_2 + 39) + (x_2^2 + 23\,x_2 + 2, y_2 + 33\,x_2 + 7) + 
    (x_2^2 + 44\,x_2 + 33, y_2 + 35\,x_2 + 28) - (x_2 + 1, y_2) - (x_2 + 11, y_2 + 42) - 
    (x_2 + 16, y_2 + 34) - (x_2 + 50, y_2 + 13) - (x_2^2 + 26\,x_2 + 12, y_2 + 49\,x_2 + 29)
    - (x_2^2 + 47\,x_2 + 5, y_2 + 7\,x + 14)$\\\hline
  \end{tabular}
\end{footnotesize}\medskip
  \caption{Some relations collected in the sieving phase.}
  \label{tab:relex}
\end{table}

\paragraph*{Linear algebra phase.}

With our smoothness choice, the factor basis is derived from places of degree
one and two.  Since we prefer functions to divisors,
 the factor basis will contain the reduction modulo
$\PP$, resp.  $\PQ$, of functions the divisors of which are equal to $76
(x_1+\alpha, y_1+\beta)-76 (1/x_1, y_1/x_1^2)$, resp. $76 (x_2+\alpha,
y_2+\beta)-76 (1/x_2, y_2/x_2^2)$ (remember that in our example $\#E(\Fp) = 76$).
In this setting, the evaluation at $\PP$ or $\PQ$ of any smooth function can
be easily written as a product of elements of the factor basis.

It is worth  recalling  that the action of the Frobenius $\phi$ on the reduction
of a function modulo $\PP$ or $\PQ$ is equal to the reduction of a function,
the poles and the zeros of which are translated by one specific point of $\Ker
I$.  In our example, this point is $F_1=(11: 48: 1)$ for the reduction
modulo $\PP$ and $F_2=(45: 34: 1)$ for the reduction modulo $\PQ$.  For
instance, let us consider a function $g_0$ the divisor of which is equal to
$76(x_1 + 41, y_1 + 8)-76(1/x_1, y_1/x_1^2)$.  Let us now consider a
function $g_6$ which corresponds to $(-41:-8~:1)+6 F_1$, that is a function
with divisor equal to $76(x_1 + 45, y_1 + 17)-76(1/x_1, y_1/x_1^2)$. We have
then
\begin{math}
  \overline{g_6} = c . 
  \overline{g_0}^{p^6}\overline{f}^{1+p+p^2+p^3+p^4+p^5}\text{for some } c\in \Fp,
\end{math}
where $f$ is a function the divisor of which is equal to $76\,F_1-76(1/x_1,
y_1/x_1^2)$.

Thanks to this observation, we can thus divide by $19$ the size of the factor
basis, at the expense in the linear algebra phase of entries equal to sums of
powers of $p$. We finally have $4$ meaningful places of degree 1 and $92$
meaningful places of degree 2 on  each side, that is a total of
$196$ entries in our factor basis. Of course, under the Galois conjugations,
most of the relations obtained in the sieving phase are redundant, but it does
not really matter since it is not difficult to reduce the sieving phase to the
only meaningful relations.

We have

$$61^{19}-1=2^2\cdot3\cdot5\cdot229\cdot607127818287731321660577427051.$$

We 
performed the linear algebra modulo the largest factor of $61^{19}-1$, that is
the $99$-bit integer $607127818287731321660577427051$.
This gives us the discrete logarithm in basis $f\bmod \cI$ of any
element in the smoothness basis.
For instance, if  $g$ is  any function such that div $g$ = $76\,(x_1^2 + 37\,x_1
+ 54, y_1+ 41\,x_1 + 16)-152(1/x_1, y_1/x_1^2)$, we find that
\begin{displaymath}
  g^{2^2\cdot3\cdot5\cdot229} = 
  (f^{2^2\cdot3\cdot5\cdot229})^{471821537021905592692223848756}.\end{displaymath}

\section{Generalization and limitations}\label{section:conclusion}

The construction in section~\ref{section:carre} can and should be
generalized. 

Let  
$E$ be again  an ordinary elliptic curve
over $\Fp$  and let  $\igot$ be an invertible ideal 
in the endomorphism ring 
 $\End(E)$. 
We assume that $\igot$ divides   $\phi-1$  and  $\End(E)/\igot$ is cyclic
of order  $d\ge 2$. Let  $F$ be the  quotient of  $E$ by the kernel $T$ of
$\igot$ and  $I: E \rightarrow F$ the quotient isogeny.

The integer  $d$ belongs to the ideal  $\igot$. Let
 $u$ and $v$ be two elements in  $\igot$ such that
$d=u+v$ and  $(u)=\igot\agot_1\bgot_1$
and  $(v)=\igot\agot_2\bgot_2$ where  $\agot_1$, $\bgot_1$, $\agot_2$, $\bgot_2$
are invertible ideals in  $\End(E)$.
We deduce the existence 
of two elliptic curves  $E_1$ and  $E_2$ and four isogenies
 $\alpha_1$, $\beta_1$, $\alpha_2$, $\beta_2$, such that 
$\beta_1\alpha_1+\beta_2\alpha_2=I$. 

We represent all these
isogenies on the (non commutative) diagram below.
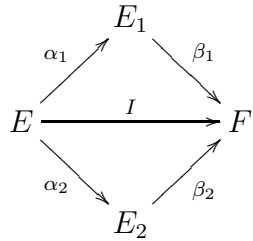
\begin{eqnarray*}\xymatrix{
 & E_1 \ar@{->}^{\beta_1}[dr]   &   \\
E \ar@{->}^{I}[rr]\ar@{->}^{\alpha_1}[ur]\ar@{->}_{\alpha_2}[dr]&  &F\\
& E_2  \ar@{->}_{\beta_2}[ur]&
}
\end{eqnarray*}

We set $\cS=E_1\times E_2$. As for  $\cA$ we choose the image of
  $(\alpha_1,\alpha_2): E\rightarrow \cS$. And 
$\cB$ is the inverse image of  $f$ by  $\beta_1+
\beta_2: \cS\rightarrow F$ where  $f$  generates
the quotient  $F(\Fp)/I(E(\Fp))$.
The intersection of  $\cA$  and  $\cB$ is the image
by  $(\alpha_1, \alpha_2)$  of  $I^{-1}(f)\subset E$.
We choose  $u$ and  $v$ such that  $\agot_1$, $\bgot_1$,
$\agot_2$,  and  $\bgot_2$, have norms close to
the square root of $d$.

This construction is useful
when the norm of $\igot$ is
much smaller than the norm of $\phi-1$.
So we managed to construct Galois invariant smoothness
basis for a range of finite fields. Our constructions go beyond
the classical Kummer case. They are efficient
when the degree $d$ is either below   $4\sqrt q$ or in the
interval 
$]q+1-2\sqrt q,q+1+2\sqrt q [$.

\bibliography{angfri7}
\end{document}